\documentclass[12pt]{article}

\usepackage{amsfonts,amssymb,amsmath,amsthm,epsfig,euscript}

\newtheorem{thm}{Theorem}[section]
\newtheorem{prop}[thm]{Proposition}
\newtheorem{cor}[thm]{Corollary}
\newtheorem{lem}[thm]{Lemma}
\newtheorem{conj}[thm]{Conjecture}
\newtheorem{exa}[thm]{Example}

\long\def\symbolfootnote[#1]#2{\begingroup
\def\thefootnote{\fnsymbol{footnote}}\footnote[#1]{#2}\endgroup}

\newcommand{\ben}{\begin{enumerate}}
\newcommand{\een}{\end{enumerate}}
\newcommand{\ble}{\begin{lem}}
\newcommand{\ele}{\end{lem}}
\newcommand{\bth}{\begin{thm}}
\renewcommand{\eth}{\end{thm}}
\newcommand{\bpr}{\begin{prop}}
\newcommand{\epr}{\end{prop}}
\newcommand{\bco}{\begin{cor}}
\newcommand{\eco}{\end{cor}}
\newcommand{\bcon}{\begin{conj}}
\newcommand{\econ}{\end{conj}}
\newcommand{\bde}{\begin{defn}}
\newcommand{\ede}{\end{defn}}
\newcommand{\bex}{\begin{exa}}
\newcommand{\eex}{\end{exa}}
\newcommand{\barr}{\begin{array}}
\newcommand{\earr}{\end{array}}
\newcommand{\btab}{\begin{tabular}}
\newcommand{\etab}{\end{tabular}}
\newcommand{\beq}{\begin{equation}}
\newcommand{\eeq}{\end{equation}}
\newcommand{\bea}{\begin{eqnarray*}}
\newcommand{\eea}{\end{eqnarray*}}
\newcommand{\bce}{\begin{center}}
\newcommand{\ece}{\end{center}}
\newcommand{\bpi}{\begin{picture}}
\newcommand{\epi}{\end{picture}}
\newcommand{\bpp}{\begin{picture}}
\newcommand{\epp}{\end{picture}}
\newcommand{\bfi}{\begin{figure} \begin{center}}
\newcommand{\efi}{\end{center} \end{figure}}
\newcommand{\bprf}{\begin{proof}}
\newcommand{\eprf}{\end{proof}\medskip}

\newcommand{\bsl}{\begin{slide}{}}
\newcommand{\esl}{\end{slide}}
\newcommand{\bfr}{\begin{frame}}
\newcommand{\efr}{\end{frame}}

\newcommand{\pf}{{\it Proof.}\quad}
\newcommand{\qmq}[1]{\quad\mbox{#1}\quad}
\newcommand{\spn}[1]{\langle{#1}\rangle}
\newcommand{\hqed}{\hfill \qed}
\newcommand{\hqedm}{\hfill \qed \medskip}

\newcommand{\ree}[1]{(\ref{#1})}
\newcommand{\emp}{\emptyset}

\newcommand{\sbe}{\subseteq}

\newcommand{\ol}{\overline}
\newcommand{\ra}{\rightarrow}

\newcommand{\fS}{{\mathfrak S}}
\newcommand{\al}{\alpha}
\newcommand{\be}{\beta}

\newcommand{\de}{\delta}
\newcommand{\ep}{\epsilon}

\newcommand{\io}{\iota}

\newcommand{\om}{\omega}

\newcommand{\si}{\sigma}

\newcommand{\ze}{\zeta}

\newcommand{\De}{\Delta}

\newcommand{\Om}{\Omega}
\newcommand{\Si}{\Sigma}

\newcommand{\bx}{{\bf x}}

\newcommand{\bbP}{{\mathbb P}}

\newcommand{\bbZ}{{\mathbb Z}}
\newcommand{\cA}{{\cal A}}

\newcommand{\cF}{{\cal F}}

\newcommand{\cL}{{\cal L}}

\newcommand{\cM}{{\cal M}}

\newcommand{\cS}{{\cal S}}

\newcommand{\wb}{\bar{w}}

\newcommand{\case}[4]{\left\{\barr{ll}#1&\mbox{#2}\\#3&\mbox{#4}\earr\right.}
\newcommand{\dil}{\displaystyle}
\newcommand{\drac}[2]{\dil\frac{#1}{#2}}
\newcommand{\wt}{\operatorname{wt}}
\newcommand{\Em}{\operatorname{Em}}
\newcommand{\adj}{\operatorname{adj}}

\newcommand{\gfo}{\operatorname{gfo}}
\newcommand{\fo}{\operatorname{fo}}
\newcommand{\clp}{\operatorname{clp}}

\newcommand{\fig}[3]
{\begin{figure}[ht] \centerline{\scalebox{#1}{\epsfig{file=#2.eps}}}
\vspace{-1mm} \caption{#3} \label{fig:#2}
\end{figure}}

\title{Rationality, irrationality, and Wilf equivalence in generalized
  factor order
}

\author{
Sergey Kitaev \\
\small Institute of Mathematics \\[-0.8ex]
\small Reykjav\'{i}k University \\[-0.8ex]
\small IS-103 Reykjav\'{i}k, Iceland\\[-0.8ex]
\small \texttt{sergey@ru.is}
\and
Jeffrey Liese \\
\small Department of Mathematics\\[-0.8ex]
\small University of California, San Diego\\[-0.8ex]
\small La Jolla, CA 92093-0112. USA\\[-0.8ex]
\small \texttt{jliese@math.ucsd.edu} \and
Jeffrey Remmel\footnote{Partially supported by NSF grant DMS 0654060} \\
\small Department of Mathematics\\[-0.8ex]
\small University of California, San Diego\\[-0.8ex]
\small La Jolla, CA 92093-0112. USA\\[-0.8ex]
\small \texttt{remmel@math.ucsd.edu}
\and
Bruce E. Sagan \\
\small Department of Mathematics\\[-0.8ex]
\small Michigan State University\\[-0.8ex]
\small East Lansing, MI 48824-1027\\[-0.8ex]
\small \texttt{sagan@math.msu.edu} }

\date{
{\small MR Subject Classifications: 05A15, 68R15, 06A07}
{\it Dedication.\/}  This paper is dedicated to Anders Bj\"orner on
the occasion of his 60th
birthday.  His work has very heavily influenced ours.
}

\begin{document}
\maketitle

\begin{abstract}

Let $P$ be a partially ordered set and consider the free monoid $P^*$
of all words over $P$.  If $w,w'\in P^*$ then $w'$ is a factor of $w$
if there are words $u,v$ with $w=uw'v$.  Define generalized factor
order on $P^*$ by letting $u\le w$ if there is a factor $w'$ of $w$
having the same length as $u$ such that $u\le w'$, where the
comparison of $u$ and $w'$ is done componentwise using the partial
order in $P$. 
One obtains ordinary factor order by insisting that $u=w'$ or,
equivalently, by taking $P$ to be an antichain.

Given $u\in P^*$, we prove that the language 
$\cF(u)=\{w\ :\ w\ge u\}$ is accepted by a finite state
automaton. If $P$ is finite then it 
follows that the generating function $F(u)=\sum_{w\ge u} w$ is 
rational.  This is an analogue of a theorem of Bj\"orner and Sagan for
generalized subword order.  

We also consider  $P=\bbP$,
the positive integers with the usual total order, so that $P^*$ is the
set of compositions.  In this case one
obtains a weight generating function $F(u;t,x)$ by substituting $tx^n$ each time
$n\in\bbP$ appears in $F(u)$.  We show that this generating
function is also rational by using the transfer-matrix method.  Words
$u,v$ are said to be Wilf equivalent 
if $F(u;t,x)=F(v;t,x)$ and we prove various Wilf equivalences
combinatorially.

Bj\"orner found a recursive formula for the M\"obius function of
ordinary factor order on
$P^*$.  It follows that one always has $\mu(u,w)=0,\pm1$.
Using the Pumping Lemma we show that
the generating function $M(u)=\sum_{w\ge u} |\mu(u,w)| w$ can be irrational.
\end{abstract}

\section{Introduction and definitions}

Let $P$ be a set and consider the corresponding {\it free monoid\/} or
{\it Kleene closure\/} of all words over $P$:
$$
P^*=\{w=w_1w_2\ldots w_\ell\ :\ \mbox{$n\ge0$ and $w_i\in P$ for all $i$}\}.
$$
Let $\ep$ be the empty word and for any $w\in P^*$ we denote its
cardinality or {\it length\/} by $|w|$.
Given $w,w'\in P^*$, we say that $w'$ is a {\it factor\/} of $w$ if
there are words $u,v$ with $w=uw'v$, where adjacency denotes
concatenation.  
For example, $w'=322$ is a factor of $w=12213221$ starting with the fifth
element of $w$.
{\it Factor order\/} on $P^*$ is the partial order
obtained by letting $u\le_{\fo} w$ if and only if  there is a factor $w'$ of
$w$ with $u=w'$.

Now suppose that we have a poset $(P,\le)$.  We define {\it
  generalized factor order\/} on $P^*$ by letting $u\le_{\gfo} w$ if there
  is a factor $w'$ of $w$ such that
\ben
\item[(a)] $|u|=|w'|$, and
\item[(b)] $u_i\le w'_i$ for $1\le i \le |u|$.
\een
We call $w'$ an {\it embedding\/} of $u$ into $w$, and
if the first element of $w'$ is the $j$th element of $w$, we call $j$
an {\it embedding index\/} of $u$ into $w$.  We also say that in this
embedding $u_i$ is in {\it position\/} $j+i-1$.
To illustrate, suppose $P=\bbP$, the positive integers with the usual
order relation.  If $u=322$ and $w=12213431$ then $u\le_{\gfo} w$ because
of the embedding factor $w'=343$ which has embedding index 5, and the
two 2's of $u$ are in positions 6 and 7.
Note that we obtain ordinary factor order by taking $P$ to be an
antichain.  Also, we will henceforth drop the subscript $\gfo$ since
context will make it clear what order relation is meant.  Generalized
factor order is the focus of this paper.

Returning to the case where $P$ is an arbitrary set,
let $\bbZ\spn{\spn{P}}$
be the algebra of formal power series with integer coefficients and
having the elements of $P$ as noncommuting variables.  In other words,
$$
\bbZ\spn{\spn{P}}=\left\{ f =\sum_{w\in P^*} c(w) w\ :\ 
\mbox{$c(w)\in\bbZ$ for all $w$}\right\}.
$$
If $f\in\bbZ\spn{\spn{P}}$ has no constant term, i.e., $c_\ep=0$, then
define
$$
f^*=\ep+f+f^2+f^3+\cdots=(\ep-f)^{-1}.
$$
(We need the restriction on $f$ to make sure that the sums are well
defined as formal power series.)   We say that $f$ is {\it rational\/}
if it can be constructed from the elements of $P$ using only a finite
number of applications of the algebra operations and the star
operation.

A {\it language\/} is any $\cL\sbe P^*$.  It has an associated
generating function 
$$
f_\cL=\sum_{w\in \cL} w.
$$
The language $\cL$ is {\it regular\/} if $f_\cL$ is rational.

Consider generalized factor order on $P^*$ and fix a word $u\in P^*$.
There is a 
corresponding language and generating function
$$
\cF(u)=\{w\ :\ w\ge u\}\qmq{and}
F(u)=\sum_{w\ge u} w.
$$
One of our main results is as follows.
\bth
\label{F(u)}
If $P$ is a finite poset and $u\in P^*$ then $F(u)$ is rational.
\eth

This is an analogue of a result of Bj\"orner and Sagan~\cite{bs:rmf}
for generalized subword order on $P^*$.  
For related results, also see Goyt~\cite{goy:mfr}.
{\it Generalized subword  order\/} is defined exactly like
generalized factor order except that $w'$ is only required to be a
subword of $w$, i.e., the elements of $w'$ need not be consecutive in $w$.

To prove the previous theorem, we will use finite automata.
Given any set, $P$, a {\it nondeterministic finite automaton\/}  or 
{\it NFA\/} over $P$ is a digraph (directed graph) $\De$ with vertices
$V$ and arcs $\vec{E}$ having the following properties.
\ben
\item  The elements of $V$ are called {\it states\/} and $|V|$ is finite.
\item There is a designated {\it initial state\/} $\al$ and a set
  $\Om$ of {\it final states}.
\item  Each arc of $\vec{E}$ is labeled with an element of $P$.
\een
Given a (directed) path in $\De$ starting at
  $\al$, we construct a word in $P^*$ by 
concatenating the elements on the arcs on the path in the order in
which they 
are encountered.  The {\it language accepted by $\De$\/} is the set
of all such words which are associated with  paths ending in a final
state.  It is a well-known theorem that, for $|P|$ finite,  a
  language $\cL\sbe P^*$ is regular if and only if there is a NFA
  accepting $\cL$.  (See, for example,  the text of
  Hopcroft and Ullman~\cite[Chapter 2]{hu:iat}.)

We will demonstrate Theorem~\ref{F(u)}
by constructing a NFA accepting the
language for $F(u)$.  This will be done in the next section.
In fact, the NFA still exists even if $P$ is infinite,
suggesting that more can be said about the generating function in this
case.  

We are particularly interested in the case of $P=\bbP$ with
the usual order relation.   So $\bbP^*$ is just the set of {\it compositions\/}
(ordered integer partitions).  Given $w=w_1 w_2\ldots w_\ell\in\bbP^*$, we
define its {\it norm\/} to be
$$
\Si(w)=w_1+w_2+\cdots+w_\ell.
$$
Let $t,x$ be commuting variables.  Replacing each $n\in w$ by $tx^n$
we get an associated monomial called the {\it weight\/} of $w$
$$
\wt(w)=t^{|w|} x^{\Si(w)}.
$$ 
For example, if $w=213221$ then
$$
\wt(w) = tx^2\cdot tx\cdot tx^3\cdot tx^2\cdot tx^2\cdot tx
= t^6 x^{11}.
$$
We also have the associated
{\it weight generating function\/}
$$
F(u;t,x)=\sum_{w\ge u} \wt(w).
$$
Our NFA will
demonstrate, via the transfer-matrix method, that this is also a
rational function of $t$ and $x$. 
The details will be given in Section~\ref{pi}.

Call $u,w\in\bbP^*$ {\it Wilf equivalent\/} if $F(u;t,x)=F(v;t,x)$.
This definition is modelled on the one used in the theory of pattern
avoidance.  See the survey article of Wilf~\cite{wil:pp} for more
information about this subject.  Section~\ref{we}  is devoted to
proving various Wilf equivalences.  Although these results were
discovered by having a computer construct the corresponding generating
functions, 
the proofs we give are purely combinatorial.  In the next two
sections, we investigate a stronger notion of equivalence and compute
generating functions for two families of compositions.

Bj\"orner~\cite{bjo:mff} gave a recursive formula for the M\"obius
function of (ordinary) factor order.   It follows from his theorem
that $\mu(u,w)=0,\pm1$ for all $u,w\in P^*$.  Using the Pumping
Lemma~\cite[Lemma 3.1]{hu:iat} we show that there are finite sets $P$
and $u\in P^*$ such that the language 
$$
\cM(u)=\{w\  :\ \mu(u,w)\neq 0\}
$$
is not regular.  This is done in  Section~\ref{mf}.  The penultimate
section is devoted to comments, conjectures, and open questions.  And
the final one contains tables.

\section{Construction of automata}
\label{ca}

We will now introduce two other languages which are related to
$\cF(u)$ and which will be useful in proving Theorem~\ref{F(u)}
and in demonstrating Wilf equivalence.  We say that $u$ is a 
{\it suffix\/} (respectively, {\it prefix\/}) of $w$ if $w=vu$
(respectively, $w=uv$) for some word $v$.  Let $\cS(u)$ be all the
$w\in \cF(u)$ such that, in the definition of generalized factor order,
the only possible choice for $w'$ is
a suffix of $w$.  Let $S(u)$ be the corresponding generating function.

We say that $w\in P^*$ {\it avoids\/}  $u$ if $w\not\geq u$ in generalized
factor order.  Let $\cA(u)$ be the associated language 
with generating function $A(u)$.  The next result follows easily from
the definitions and so we omit the proof.  In it, we will use the notation
$Q$ to stand both for a subset of $P$ and for the generating function
$Q=\sum_{a\in Q} a$.  Context will make it clear which is meant.
\ble
\label{cFSA}
Let $P$ be any poset and let $u\in P^*$.  Then we have the following
relationships:  
\ben
\item $\cF(u)=\cS(u) P^*$ and $F(u)=S(u) (\ep-P)^{-1}$,
\item $\cA(u)=P^*-\cF(u)$ and $A(u)=(\ep-P)^{-1} - F(u)$.\hqed
\een
\ele

We will now prove that all three of the languages we have defined are
accepted by NFAs.  An example follows the proof so the reader may want
to read it in parallel.
\bth
\label{NFA}
Let $P$ be any poset and let $u\in P^*$.  Then there are NFAs
accepting $\cF(u)$, $\cS(u)$, and $\cA(u)$.
\eth
\pf
We first construct an NFA, $\De$, for $\cS(u)$.  Let $\ell=|u|$.  The
states of $\De$ will be all subsets $T$ of $\{1,\ldots,\ell\}$.  The
initial state is $\emp$.  Let 
$w=w_1\ldots w_m$ be the word corresponding to a path from $\emp$ to $T$.
The NFA will be constructed so that if the path is continued, the 
only possible  embedding indices are those in
the set $\{m-t+1\ :\ t\in T\}$.  In other words, for each $t\in T$ we
have 
\beq
\label{uw}
u_1 u_2\ldots u_t \le w_{m-t+1} w_{m-t+2} \ldots w_m,
\eeq
for each $t\in\{1,\ldots,\ell\}-T$ this inequality does not hold,
and $u\not\le w'$ for any factor $w'$ of $w$ starting at an index smaller
then $m-\ell+1$.  From this description it is clear that the final
states should be those containing  $\ell$.

The definition of the  arcs of $\De$ is forced by the interpretation
of the states.  There will be no arcs out of a final state.
If $T$ is a nonfinal state and $a\in P$ then there will be an arc from
$T$ to 
$$
T'=\{t+1\ :\ \mbox{$t\in T\cup\{0\}$ and $u_{t+1}\le a$}\}.
$$
It is easy to see that~\ree{uw} continues to hold for all $t'\in T'$
once we append $a$ to $w$.  This finishes the construction of the NFA
for $\cS(u)$.

To obtain an automaton for $\cF(u)$, just add loops to the final
states of $\De$, one for each $a\in P$.  An automaton for $\cA(u)$ is
obtained by just interchanging the final and nonfinal states in the
automaton for $\cF(u)$.\hqedm

As an example, consider $P=\bbP$ and $u=132$.   We will do several
things to simplify writing down the automaton.  First of all, certain
states may not be reachable by a path starting at the initial state.
So we will not display such states.  For example, we can not reach the
state $\{2,3\}$ since $u_1=1\le w_i$ for any $i$ and so $1$ will be in
any state reachable from $\phi$.
Also, given states $T$ and $U$
there may be many arcs from $T$ to $U$, each having a different label.  So
we will replace them by one arc bearing the set of labels of all such
arcs.  Finally, set braces will be dropped for readability.  The
resulting digraph is displayed in Figure~1. 

\fig{.7}{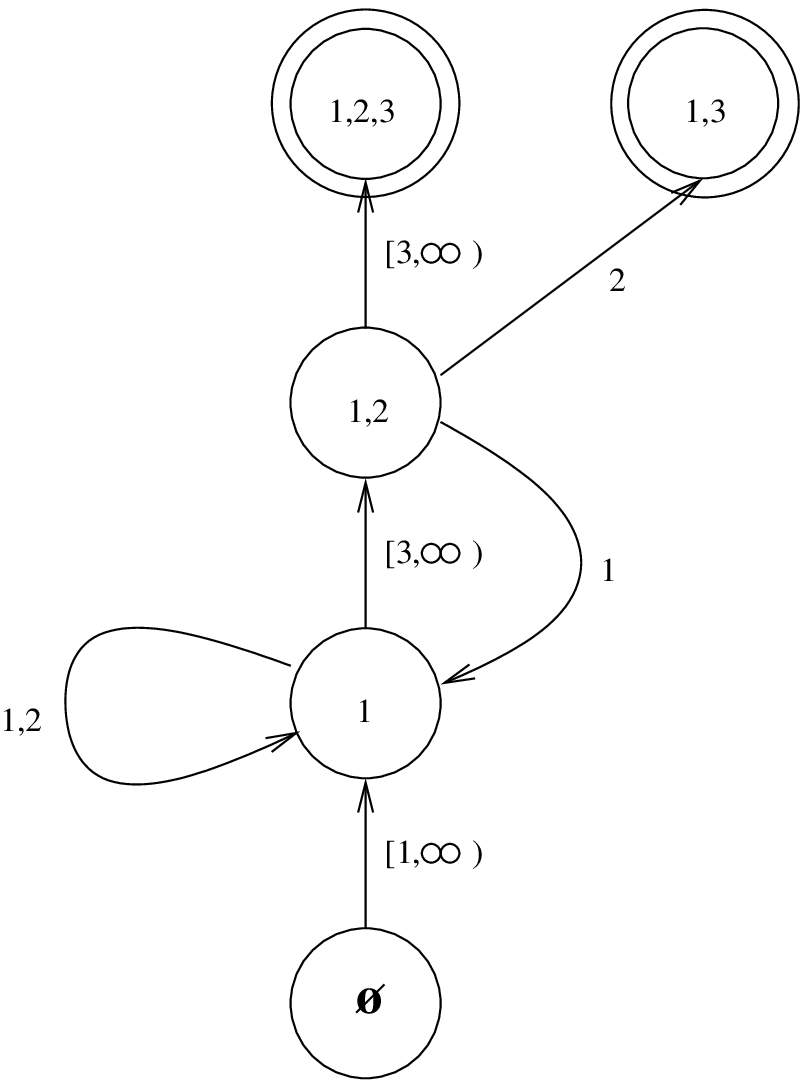}{A NFA accepting $\cS(132)$}

Consider what happens as we build a word $w$ starting from the initial
state $\emp$.
Since $u_1=1$, any element of $\bbP$ could be the first element
of an embedding of $u$ into $w$.  That is why every element of the interval
$[1,\infty)=\bbP$ produces an arrow from the initial state to the
state $\{1\}$.  Now if $w_2\le 2$, then an embedding of $u$ could no
longer start at $w_1$ and so these elements give loops at the state
$\{1\}$. But if $w_2\ge3$ then an embedding could start at either
$w_1$ or at $w_2$ and so the corresponding arcs all go to the state
$\{1,2\}$.  The rest of the automaton is explained similarly.

As an immediate consequence of the previous theorem we get the
following result which includes Theorem~\ref{F(u)}.
\bth
\label{FSA}
Let $P$ be a finite poset and let $u\in P^*$.  Then the generating
functions $F(u)$, $S(u)$, and $A(u)$ are all rational.\hqed
\eth

\section{The positive integers}
\label{pi}

If $P=\bbP$ then Theorem~\ref{FSA} no longer applies to the generating
functions $F(u)$, $S(u)$, and $A(u)$.  However, we can still show
rationality of the weight generating function $F(u;t,x)$ as defined
in the introduction.  Similarly, we will see that the series
$$
S(u;t,x)=\sum_{w\in\cS(u)} \wt(w)\qmq{and}
A(u;t,x)=\sum_{w\in\cA(u)} \wt(w)
$$
are rational.

Note first that Lemma~\ref{cFSA} still holds for $\bbP$ and can be
made more explicit in this case.  Extend the function $\wt$ to all of
$\bbZ\spn{\spn{\bbP}}$ by letting it act linearly.  Then
\bea
\wt(\ep-\bbP)^{-1}
&=& \frac{1}{1-\sum_{n\ge 1} tx^n}\\[5pt]
&=& \frac{1}{1-tx/(1-x)}\\[5pt]
&=& \frac{1-x}{1-x-tx}.
\eea
We now plug this into the lemma.
\bco
\label{FSAtx}
We have
\smallskip
\ben
\item $\dil F(u;t,x)=\frac{(1-x)S(u;t,x)}{1-x-tx}$ and
\smallskip
\item $\dil A(u;t,x)=\frac{1-x}{1-x-tx}-F(u;t,x)$.\hqed
\een
\eco
\medskip
It follows that if any one of these three series is rational then the
other two are as well.

We will now use the NFA, $\De$, constructed in Theorem~\ref{NFA} to show that
$S(u;t,x)$ is rational.  This is essentially an application of the
transfer-matrix method.  See the text of 
Stanley~\cite[Section 4.7]{sta:ec1} for more information about this
technique.  The {\it transfer matrix\/} $M$ for $\De$ has rows and columns
indexed by the states with
$$
M_{T,U}=\sum_n \wt(n)
$$
where the sum is over all $n$ which appear as labels on the arcs from
$T$ to $U$.  For example, consider the case  where $w=132$ as done at
the end of the previous section.  If we list the states in the order
$$
\emp,\ \{1\},\ \{1,2\},\ \{1,3\},\ \{1,2,3\}
$$
then the transfer matrix is
$$
M=
\left[
\barr{ccccc}
0     &\drac{tx}{1-x}&0               &0   &0\\
0     &t(x+x^2)      &\drac{tx^3}{1-x}&0   &0\\
0     &tx            &0               &tx^2&\drac{tx^3}{1-x}\\[8pt]
0     &0             &0               &0   &0\\[8pt]
0     &0             &0               &0   &0\\
\earr
\right]
$$

Now $M^k$ has entries $M^k_{T,U}=\sum_{w} \wt(w)$ where the sum is
over all words $w$ corresponding to a directed walk of length $k$ from
$T$ to $U$.  So to get the weight generating function for walks of
all lengths one considers $\sum_{k\ge0} M^k$.  Note that this sum
converges in the algebra of matrices over the formal power series
algebra $\bbZ[[t,x]]$ because none of the entries of $M$ has a
constant term.  It follows that
\beq
\label{LM}
L:=\sum_{k\ge0} M^k=(I-M)^{-1}=\frac{\adj(I- M)}{\det (I-M)}
\eeq
where $\adj$ denotes the adjoint.  

Now 
$$
S(u;t,x)=\sum_{T} L_{\emp,T}
$$
where the sum is over all final states of $\De$.  So it suffices to
show that each entry of $L$ is rational.  From equation~\ree{LM}, this
reduces to showing that each entry of $M$ is rational.  So consider
two given states $T,U$.  If $T$ is final then we are done since the
$T$th row of $M$ is all zeros.  If $T$ is not final, then consider
\beq
\label{T'}
T'=\{t+1\ :\ t\in T\cup\{0\}\}.
\eeq
If $U=T'$ then there will be an $N\in\bbP$ such that all the arcs out
of $T$ with labels $n\ge N$ go to $T'$.  So $M_{T,T'}$ will contain
$\sum_{n\ge N} tx^n= tx^N/(1-x)$ plus a finite number of other terms
of the form $tx^m$.  Thus this entry is rational.  If $U\neq T'$, then
there will only be a finite number of arcs from $T$ to $U$ and so
$M_{T,U}$ will actually be a polynomial.  This shows that every entry
of $M$ is rational and we have proved, with the aid of the remark
following Corollary~\ref{FSAtx}, the following result.
\bth
\label{FSAtxthm}
If $u\in\bbP^*$ then $F(u;t,x)$, $S(u;t,x)$, and $A(u;t,x)$ are all
rational.\hqed
\eth

\section{Wilf equivalence}
\label{we}

Recall that $u,v\in\bbP^*$ are {\it Wilf equivalent\/}, written 
$u\sim v$, if $F(u;t,x)=F(v;t,x)$.  By Corollary~\ref{FSAtx}, this
is equivalent to $S(u;t,x)=S(v;t,x)$ and to $A(u;t,x)=A(v;t,x)$.  It
follows that to prove Wilf equivalence, it suffices to find a
weight-preserving bijection $f:\cL(u)\ra\cL(v)$ where $\cL=\cF$,
$\cS$, or $\cA$.  Since $\sim$ is an equivalence relation, we can
talk about the {\it Wilf equivalence class} of $u$ which is
$\{w\ :\ w\sim u\}$.
It is worth noting that the automata for the words in a Wilf
equivalence class need not bear a resemblance to each other.

Part of the motivation for this section is to try to explain as many
Wilf equivalences as possible between permutations.  For reference, in
Section~\ref{t} the first table lists all such equivalences up
through 5 elements.

First of all, we consider three operations on words in $\bbP^*$.  The 
{\it  reversal\/} of $u=u_1\ldots u_\ell$ is $u^r=u_\ell\ldots u_1$.  
It will also be of interest to consider 
$1u$, the word gotten by prepending one to $u$.
Finally, we will look at $u^+$ which is gotten by
increasing each element of $u$ by one, as well as $u^-$ which performs
the inverse operation whenever it is defined.  For some of our proofs,
it will also be useful to have the following factorization.  Given
$k\in\bbP$ and $w\in \bbP^*$ the {\it $k$-factorization\/} of $w$ is
the unique expression
$$
w=y_1\ z_1\ y_2\ z_2\ \ldots\ z_{m-1}\ y_m
$$
where $y_i\in[1,k)^*$ and $z_i\in[k,\infty)^*$ for all $i$, and all
factors are nonempty with the possible exception of $y_1$ and $y_m$.
\ble
\label{sim}
We have the following Wilf equivalences.
\ben
\item[(a)]  $u\sim u^r$,
\item[(b)]  if $u\sim v$ then $1u\sim 1v$,
\item[(c)]  if $u\sim v$ then $u^+\sim v^+$.
\een
\ele
\pf
(a)  It is easy to see that the map $w\mapsto w^r$ is a
weight-preserving bijection $\cF(u)\ra\cF(u^r)$.

(b)  We can assume we are given a weight-preserving bijection
$f:\cS(u)\ra\cS(v)$. 
Since 1 is the minimal element of $\bbP$,
$$
\cS(1u)=\{w\in\cS(u)\ :\ |w|>|u|\}.
$$
So $f$ restricts to a weight-preserving bijection from
$\cS(1u)$ to $\cS(1v)$.

(c)  Now we consider a weight-preserving bijection 
$g:\cA(u)\ra\cA(v)$.  Given $w\in\bbP^*$, let
$$
w=y_1\ z_1\ y_2\ z_2\ \ldots\ z_{m-1}\ y_m
$$ 
be its 2-factorization.
Since all elements of $u^+$ are at least two, $w\in\cA(u^+)$ if
and only if $z_i\in\cA(u^+)$ for all $i$.  This is equivalent to
$z_i^-\in\cA(u)$ for all $i$.  Thus if we map
$w$ to
$$
y_1\ g(z_1^-)^+\ y_2\ g(z_2^-)^+\ \ldots\ g(z_{m-1}^-)^+\ y_m
$$
then we will get the desired weight-preserving bijection
$\cA(u^+)\ra\cA(v^+)$.\hqedm

We can combine these three operations to prove more complicated Wilf
equivalences.  Since a word $w\in\bbP^*$ is just a sequence of
positive integers, terms like ``weakly increasing'' and ``maximum''
have their usual meanings.  Also, let $w^{+m}$ be the result of applying
the + operator $m$ times.   By using the previous lemma and
induction, we obtain the following result.  The proof is so
straight forward that it is omitted. 
\bco
\label{id}
Let $y,y'$ be weakly increasing compositions and $z,z'$ be weakly decreasing
compositions such that $yz$ is a rearrangement of $y'z'$.
Then for any $u\sim v$  we have
$$
yu^{+m}z\sim y'v^{+m} z'
$$ 
whenever $m\ge\max\{y,z\}-1$.
\hqed
\eco

Applying the two previous results, we can obtain the  Wilf equivalences in the
symmetric group $\fS_3$ of 
all the permutations of $\{1,2,3\}$:
$$
123\sim 321\sim 132\sim 231\qmq{and}213\sim 312.
$$
These two groups are indeed in different equivalence classes as one
can use equation~\ree{LM} to compute that
$$
S(123;t,x)=\frac{t^3 x^6}{(1-x)^2 (1-x -tx+ tx^3 -t^2x^4)}
$$
while
$$
S(213;t,x)=\frac{t^3x^6(1+tx^3)}{(1-x)(1-x+t^2x^4)(1-x-tx+tx^3-t^2x^4)}.
$$

However, we will need a new result to explain some of the equivalences
in $\fS_4$ such as $2134\sim 2143$.  Let $u$ be a composition such
that $\max u$ only occurs once.  Define a {\it pseudo-embedding\/}
of $u$ into $w$ to be a factor $w'$ of $w$ satisfying the two
conditions for an embedding except that the inequality may fail at
the position(s) of $\max u$.   In particular, embeddings are
pseudo-embeddings. 

An example of the construction used in the
next theorem follows the proof and can be read in parallel.

\bth
\label{mn}
Let $x,y,z\in\{1,\ldots,m\}^*$ and suppose $n>m$.  Then
$$
xmynz\sim xnymz.
$$
\eth
\pf
Let $u=xmynz$ and $v=xnymz$.  We will construct a
weight-preserving bijection $\cA(u)\ra\cA(v)$.  To do this, it
suffices to construct such a bijection between the set differences
$\cA(u)-\cA(v)\ra\cA(v)-\cA(u)$ since the identity map can be used on
$\cA(u)\cap\cA(v)$.  

Given $w\in\cA(u)-\cA(v)$, consider the set
$$
\eta(w)=\{i\ :\ \mbox{there is an embedding of $v$ into $w$ with the
  $n$ in position $i$}\}.
$$
For such $i$, $w_i\ge n$.  It must also be that $w_{i+k}$ is in the
interval $[m,n)$ where
$k=|y|+1$:  Certainly $w_{i+k}\ge m$ because of the embedding.  But if
$w_{i+k}\ge n$ then there would also be an embedding of $u$ at the same
position as the one for $v$, contradicting $w\in\cA(u)$.

Now for each $i\in\eta(w)$ we consider the {\it string\/} beginning at $i$
$$
\si(i)=\{i,\ i+k,\ i+2k,\ \ldots,\ i+\ell k\}
$$ 
where $\ell $ is the least nonnegative integer such that there is no
pseudo-embedding 
of $v$ into $w$ with the $n$ in  position $i+\ell k$.  Note that $\ell $
depends on $i$ even though this is not  reflected in our notation.
Also, $\ell \ge1$ since there is embedding of $v$ into $w$ with the $n$ in
position $i$.  Finally, it is easy to see that
$w_{i+k},w_{i+2k},\ldots,w_{i+\ell k}\in[m,n)$ by an argument similar to that
for $w_{i+k}$.  This implies that any two strings are disjoint since
$w_i\ge n$ for $i\in\eta(w)$.

Now map $w$ to $\wb$ which is constructed by switching the values of
$w_i$ and $w_{i+\ell k}$ for every $i\in\eta(w)$.  Since strings are
disjoint, the switchings are well defined.  We must show that
$\wb\in\cA(v)-\cA(u)$.  We prove that $\wb\in\cA(v)$ by contradiction.
The switching operation removes every embedding of $v$ in $w$.  If a
new embedding was created then, because only elements of size at least
$m$ move, the $n$ in $v$ must correspond to $\wb_{i+\ell k}$ for some
$i\in\eta(w)$.  But now there is a pseudo-embedding of $v$ into $w$ 
with the $n$ in position $i+\ell k$, contradicting the
definition of $\ell $. 

To show $\wb\not\in\cA(u)$, we will actually prove the stronger
statement that there is an embedding of $u$ in $\wb$ with the $n$ in
position $i+\ell k$ for each $i\in\eta(w)$ and these are the only
embeddings.  These embeddings exist because there is a pseudo-embedding
of $v$ into $w$ with the $n$ in position $i+(\ell -1)k$,
  $\wb_{i+\ell k}\ge n$, 
and only elements of size at least $m$ move in passing from $w$ to
$\wb$.  They are the only ones because $w\in\cA(u)$ and so any
embedding of $u$ in $\wb$ would have to have the $n$ in a position of 
the form $i+\ell k$.

Finally, we need to show that this map is bijective.  But modifying
the above construction by exchanging
the roles of $u$ and $v$ and building the strings from right to left
gives an inverse.  This completes the proof.
\hqedm

By way of illustration, suppose  $u= 1~3~5~2~4~6~3$
and $v = 1~3~6~2~4~5~3$ so that $m=5$, $n=6$, and $k=3$.  We will write our example $w$
in two line form with the upper line being the positions:
$$
w=\barr{*{22}{c}}
1&2&3&4&5&6&7&8&9&10&11&12&13&14&15&16&17&18&19&20&21&22\\
1&1&2&4&8&3&9&5&4&5 &5 &4 &5 &5 &3 &3 &3 &6 &6 &5 &5 &3.
\earr
$$
Now there are three embeddings of $v$ (and none of $u$) into $w$ with
the $6$ in positions $\eta(w)=\{5,7,18\}$.  For $i=5$ we have the
string $\si(5)=\{5,8,11,14\}$ since there are pseudo-embeddings of $v$
with the $n$ in positions $5,8,11$ but not in position $14$.
Similarly $\si(7)=\{7,10,13\}$ and $\si(18)=\{18,21\}$.  So $\wb$ is
obtained by switching $w_5$ with $w_{14}$, $w_7$ with $w_{13}$, and
$w_{18}$ with $w_{21}$ to obtain
$$
\wb=\barr{*{22}{c}}
1&2&3&4&5&6&7&8&9&10&11&12&13&14&15&16&17&18&19&20&21&22\\
1&1&2&4&5&3&5&5&4&5 &5 &4 &9 &8 &3 &3 &3 &5 &6 &5 &6 &3.
\earr
$$

It is now easy to verify that our results so far suffice to explain
all the Wilf equivalences in symmetric groups up through $\fS_4$.
They also explain most, but not all, of the ones in $\fS_5$.  We will
return to the $n=5$ case in the section on open questions.

One might wonder about the necessity of the requirement that the two
equivalent words in Theorem~\ref{mn} have a unique maximum.  However, one
can see from Table 2 in Section~\ref{t} that $122$ and $212$ are not Wilf equivalent.  So
if there is an analogue of this theorem for more general words,
another condition will have to be imposed.

One might also hope that it would be possible to do without the
strings in the proof and merely switch $w_i$ and $w_{i+k}$ for all
$i\in\eta(w)$ to get $\wb$.  This would only be invertible if the
embedding indices 
for $v$ in $w$ would be the same as those for $u$ in $\wb$.
Unfortunately, this does not always work as the following example
shows.  Consider $u=231$, $v=321$, and all $w$ which are permutations
of $1223$.  Then the members of $\cA(u)-\cA(v)$ are $1322$,
$3212$, and $3221$; while those of $\cA(v)-\cA(u)$ are $1232$, $2313$,
and $2231$.  The embedding indices of $v$ in the first three
compositions are 2, 1, and 1 (respectively); while those of $u$ in the
second three are 2, 1, and 2.  Thus preservation of the indices is not
possible in this case.  However, it would be interesting to know
when one can leave the indices invariant and this will be investigated
in the next section.

The reader may have noted that a number of the maps constructed in
proving the results of this section involve rearrangement of the letters of
the word (which makes the map automatically weight preserving).  We will now
show that if one strengthens the hypothesis of Lemma~\ref{sim} (c) by
adding a rearrangement assumption, then one can also strengthen the
conclusion by applying any strictly increasing function to $u$ and $v$.
To state and prove this result, we first need some definitions.

Say that a map $f:P^*\ra P^*$ is a {\it rearrangement\/} if
$f(w)$ is a rearrangement of $w$ for all $w\in P^*$.  Now let
$u,v\in\bbP^*$ be given.  If
$f:\bbP^*\ra\bbP^*$ is a weight-preserving bijection such that, for
all $w\in\bbP^*$, 
\beq
\label{wit}
u\le w\iff v\le f(w)
\eeq
then we say that $f$ {\it witnesses\/} the Wilf equivalence $u\sim v$.

Given any function $\io:\bbP\ra\bbP$ we extend $\io$ to $\bbP^*$ by
letting
$$
\io(u_1 u_2 \ldots u_n)=\io(u_1) \io(u_2)\ldots \io(u_n).
$$
Now assume that $\io$ is strictly increasing on $\bbP$ with range $\{k_1<k_2<\ldots\}$.
Given a word $w=w_1\ldots w_m$ in $(\bbP-[1,k_1))^*$ we form its 
{\it collapse\/}, $\clp(w)$, by replacing each letter of $w$ in the
interval $[k_j,k_{j+1})$ by $j$ for all $j\in \bbP$.  For example, if
$\io(1)=3$, $\io(2)=5$, $\io(3)=8$, and $\io(4)=13$ then
$\clp(356749438)=122213113$.  For any $u,w\in\bbP^*$, we have
\beq
\label{io}
\io(u)\le w\iff u\le\clp(w).
\eeq
We now have everything in place for proof of the next result which
resembles the proof of Lemma~\ref{sim} (c).
\bth
Suppose $u,v\in\bbP^*$ such that there is a rearrangement $f:\bbP^*\ra\bbP^*$
witnessing $u\sim v$.  Then for any strictly increasing function
$\io:\bbP\ra\bbP$ there is a rearrangement $g:\bbP^*\ra\bbP^*$
witnessing $\io(u)\sim\io(v)$.
\eth
\pf
It suffices to construct a bijective rearrangement  $g$
satisfying~\ree{wit} since then it must also be weight preserving.
Given $w\in\bbP^*$, let
$$
w = y_1\ z_1\ y_2\ z_2\ \ldots\ z_{m-1} y_m 
$$
be its $k_1$-factorization where $k_1=\io(1)$.  Clearly $\io(u)\le w$ if
and only if $\io(u)\le z_i$ for some $i$.  For each $i$, define
$$
z_i'=f(\clp(z_i)).
$$
By our assumptions and~\ree{io} we have
$$
\io(u)\le z_i\iff u\le\clp(z_i)\iff v\le z_i'.
$$

Now fix $j\ge 1$ and let $z_i(1)\ldots z_i(r_j)$ be the elements of
$z_i$ in $[k_j,k_{j+1})$, reading from left to right.  These are the
elements of $z_i$ which get replaced by $j$ when passing from $z_i$
to $\clp(z_i)$.  Since $z_i'=f(\clp(z_i))$ is a rearrangement of
$\clp(z_i)$, there must be $r_j$ occurrences of $j$ in $z_i'$.  Replace
these $j$'s by $z_i(1)\ldots z_i(r_j)$, reading from left to right.
Do this for each $j\in\bbP$ and call the result $g(z_i)$.  Then
$g(z_i)$ is a rearrangement of $z_i$ and $\clp(g(z_i))=z_i'$.  It
follows from~\ree{io} and the previous displayed equation that
$$
\io(v)\le g(z_i)\iff v\le z_i'\iff \io(u)\le z_i.
$$

Now let
$$
g(w) =  y_1\ g(z_1)\ y_2\ g(z_2)\ \ldots\ g(z_{m-1}) y_m.
$$
This map is a rearrangement by construction and satisfies~\ree{wit}
because of the last displayed equation in the previous paragraph.
One can construct $g^{-1}$ from $f^{-1}$ in the same way that we
constructed $g$ from $f$.  So we are done.\hqed

\section{Strong Wilf equivalence}
\label{swe}

Given $v,w\in\bbP^*$ we let
$$
\Em(v,w)=\{j\ :\ \mbox{$j$ is an embedding index of $v$ into $w$}\}.
$$
Call compositions $u,v$ {\it strongly Wilf equivalent\/}, written
$u\sim_s v$, if there is a weight-preserving bijection
$f:\bbP^*\ra\bbP^*$ such that 
\beq
\label{Em}
\Em(u,w)=\Em(v,f(w))
\eeq
for all $w\in\bbP^*$.  
In this case we say that $f$ {\it witnesses\/} the strong Wilf
equivalence $u\sim_s v$.  Clearly strong Wilf equivalence implies Wilf
equivalence.
In addition to being a natural notion, our
interest in this concept is motivated by the fact that we were able to
prove Theorem~\ref{rep} below only under the assumption of strong Wilf
equivalence, although we suspect it is true for ordinary Wilf
equivalence.  First, however, we will prove analogues of some of our
results from the previous section in this setting.

\ble
\label{sims}
If $u\sim_s v$ then
\ben
\item[(a)]  $1u\sim_s 1v$,
\item[(b)]  $1u\sim_s v1$, 
\item[(c)]  $u^+\sim_s v^+$.
\een
\ele
\pf
Let $f:\bbP^*\ra\bbP^*$ be a map satisfying~\ree{Em}.
Define maps $g:\bbP^*\ra\bbP^*$ and
$h:\bbP^*\ra\bbP^*$ by $g(\ep)=h(\ep)=\ep$ and,
for $w=by$ with $b\in\bbP$,
$$
g(by)=bf(y) \qmq{and} h(by)=f(y)b.
$$
It follows easily that these functions establish (a) and (b).
Finally, the construction used in the proof of (c) in Lemma~\ref{sim}
can be carried over to prove the analogous case here.  That is, if one
assumes that the function $g$ given there also satisfies~\ree{Em} then
the derived map will demonstrate that $u^+\sim_s v^+$.
\hqedm

As before, we can combine the previous result and induction to get a
more general equivalence.
\bco
\label{yz}
Let $y,y'$ be weakly increasing compositions and $z,z'$ be weakly decreasing
compositions such that $yz$ is a rearrangement of $y'z'$.
Then for any $u\sim_s v$  we have
$$
yu^{+m}z\sim_s y'v^{+m} z'
$$ 
whenever $m\ge\max\{y,z\}-1$.
\hqed
\eco

Not every Wilf equivalence is a strong Wilf equivalence.  From
Lemma~\ref{sim} (a) we know that $w\sim w^r$.  But we can show that
$2143\not\sim_s 3412$ as follows.  Consider how one could construct a
word $w$ of length 7 with $\Si(w)$ minimum and
$\Em(2143,w)=\{1,3,4\}$.  Construct a table with a copy of $2143$
starting in the first, third, and fourth positions in rows 1, 2, and 3,
respectively.  Then take the maximum value in each column for the
corresponding entries of $w$:
$$
\barr{ccccccccc}
 & &2&1&4&3\\
 & & & &2&1&4&3\\
 & & & & &2&1&4&3\\
\hline
w&=&2&1&4&3&4&4&3.
\earr
$$
By construction, $w$ has the desired embedding indices and one sees
immediately that it has no others.  Note that this is the unique $w$
satisfying the given restrictions and that $\wt(w)=t^7 x^{21}$.  But
applying the same process to $3412$ gives $\wb=3434422$ with
$\wt(\wb)=t^7 x^{22}$.  Since the weights do not agree, we can not have
strong Wilf equivalence.

Finally, we come to the result alluded to at the beginning of this section.
Given $b\in\bbP$ we let $b^k$ denote the composition consisting of $k$
copies of $b$.
\bth
\label{rep}
Suppose $u=u_1\ldots u_n \sim_s v=v_1\ldots v_n$.
Then for any $k\in\bbP$
$$
u_1^k\ldots u_n^k\sim_s v_1^k\ldots v_n^k.
$$
\eth
\pf
Let $f:\bbP^*\ra\bbP^*$ be a map satisfying~\ree{Em}.
Given any $w\in\bbP^*$ and $i$ with $1\le i\le k$, consider the
subword $w[i]=w_i w_{i+k} w_{i+2k}\ldots$ of $w$.
Then the embeddings  of $u_1^k\ldots u_n^k$ in $w$ are completely determined
by the embeddings of $u$ in the $w[i]$ and vice-versa.  So
replacing each subword $w[i]$ by the subword $f(w[i])$ yields the
desired map.
\hqedm

Just as in the previous section, we can get an interesting result by
imposing  the rearrangement condition on maps.  Here is
an analogue of Corollary~\ref{yz} in this setting without the weakly
increasing assumption.
\bth
Fix $k\in\bbP$ and suppose $u,v\in[k,\infty)^*$ such that there is a
  rearrangement $f:\bbP^*\ra\bbP^*$ 
witnessing $u\sim_s v$.  Then for any two words $y,z\in[1,k]^*$
there is a rearrangement $g:\bbP^*\ra\bbP^*$
witnessing $yuz\sim_s yvz$.
\eth
\pf
It suffices to construct a bijective rearrangement  $g$
satisfying~\ree{Em} since then it must also be weight preserving.
Given $w\in\bbP^*$, let
$$
w = \psi_1\ \om_1\ \psi_2\ \om_2\ \ldots\ \om_{m-1} \psi_m 
$$
be its $k$-factorization.
Define
$$
w'=g(w)=\psi_1\ f(\om_1)\ \psi_2\ f(\om_2)\ \ldots\ f(\om_{m-1}) \psi_m.
$$
This is clearly a bijective rearrangement, so we just need to verify~\ree{Em}.

If $yuz$ embeds in $w$ at some index,
then we must show $yvz$ embeds in $w'$ at the same index.
(Showing the converse is similar.)
Now $yuz\le w$ if and only if $u\le \om_i$ for some $i$.  By assumption,
$v$ embeds in $f(\om_i)$ at the same index.  We will show that $y$
embeds in $w'$ just before this embedding of $v$.  (The proof that $z$ embeds
just after is similar.)  So consider any element $y_p$ with 
$y_p\le w_q$ in the embedding of $yuz$ in $w$.  If $w_q\in\psi_j$ for
some $j$, then $w_q'=w_q\ge y_p$.  If $w_q\in\om_j$ for some $j$, then
$w_q'\ge k\ge y_p$ because $f$ is a rearrangement.  So $y_q$ will still
embed at index $q$ in $w'$.  Thus $yvz$ embeds in $w'$ as desired
and we have completed the proof.\hqedm

As an application of this theorem, we will derive a strong Wilf
equivalence in $\fS_5$ which we could not obtain from our previous
results alone.  The proof of Lemma~\ref{sims}~(b) shows that
$123\sim_s 231$ is witnessed by a rearrangement.  From this and the proof
of Lemma~\ref{sims}~(c), it follows that 
$345\sim_s 453$ is witnessed by a rearrangement.
So the theorem just proved shows that $34512\sim_s 45312$.

\section{Computations}
\label{c}

We will now explicitly calculate the generating functions $S(u;t,x)$
for two families of words $u$.  Aside from providing an application of
the ideas from the previous sections, these particular power series
are of interest because they have numerators which are single
monomials.  This is not always the case.  For example,
$$
S(212;t,x)=\frac{t^3x^5(1+tx^2)}{(1-x)(1-x+t^2x^3)(1-x-tx+tx^2-t^2x^3)}.
$$
One can use the theory of Gr\"obner bases to show
that $(1-x)(1-x+t^2x^3)(1-x-tx+tx^2-t^2x^3)$ is not in the ideal
generated by $1+tx^2$. So $1+tx^2$ does not divide
$(1-x)(1-x+t^2x^3)(1-x-tx+tx^2-t^2x^3)$ and we can not write
$S(212;t,x)$ in the form 
$t^ax^b/Q(x,t)$ for some polynomial $Q(x,t)$.

We first determine the generating function for increasing
permutations.  It will be convenient to have the standard notation
that, for a nonnegative integer $k$,
$$
[k]_x=1+x+x^2+\cdots+x^{k-1}.
$$
\bth
For $n\ge2$, define polynomials $B_n(t,x)$ by
\bea
B_2(t,x)&=&tx(1-x)^2,\\
B_{n+1}(t,x)&=&tx^{n+1}B_n(t,x)+tx(1-x)^n(1-x^n).
\eea
Then
$$
S(12\ldots n;t,x)=\frac{t^n x^{{n+1\choose 2}}}{(1-x)^n-B_n(t,x).}
$$
\eth
\pf
Since $12\ldots n\sim n\ldots 21$, it suffices to compute the
generating function for the latter.  In that case, one can
simplify the automaton $\De$ constructed in Theorem~\ref{NFA}.

Note that $T$ is an accepting state for $\De$ if and only if 
$\max T = n$ (where 
we define $\max \emp= 0$).  Furthermore, because of our choice of
permutation, if  there is an arc from $T$ to $U$ labeled $a$, then
$\max U$ is completely determined by $\max T$ and $a$.  So we can
contract all the states with the same maximum into one.  And when we
do so, arcs of
the same label will collapse together.  The result for $n=5$ is shown in
Figure~2.  For convenience in later indexing, the
state labeled $k$ is the one resulting from amalgamating those with
maximum $n-k$.

\fig{.7}{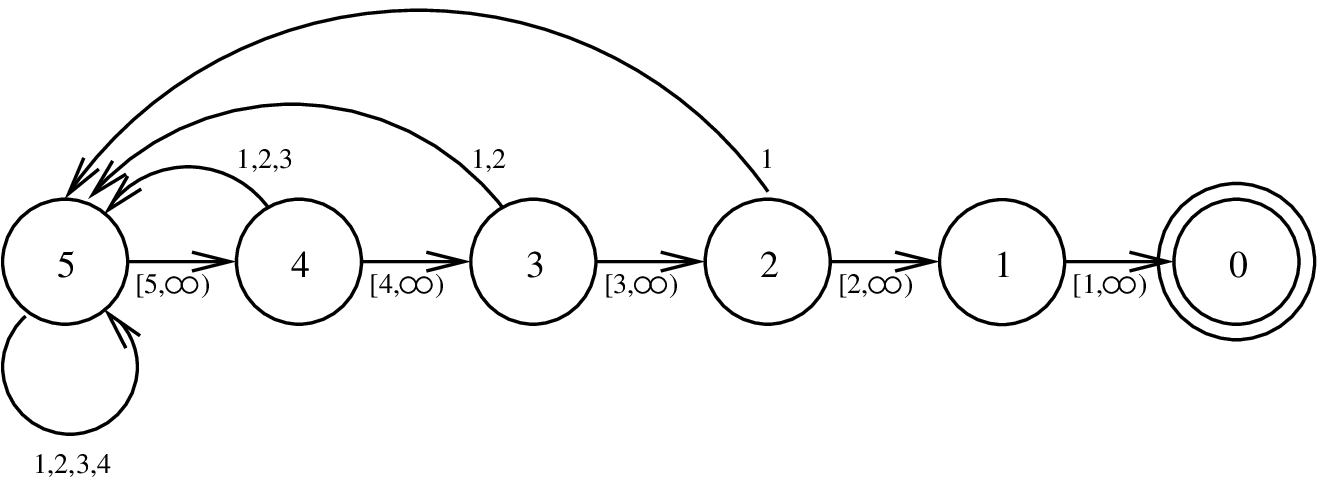}{An automaton accepting $\mathcal{S}(54321)$.}

Let $\cL_k$ be the language of all words $u$ such that the path for
$w$ starting at state $k$ leads to the accepting state $0$.  Consider
the corresponding generating function
$L_k=\sum_{u\in\cL_k} \wt(u)$.  Directly from the automaton, we have
$L_0=1$ and
$$
L_k = \frac{tx^k}{(1-x)}L_{k-1} + tx[k-1]_x L_n
$$
for $k\ge1$.  It is now easy to prove by induction that, for $k\ge2$,
$$
L_k=\frac{t^k x^{{k+1\choose2}} + B_k(t,x) L_n}{(1-x)^k}.
$$
Plugging in $k=n$ and solving for $L_n=S(n\ldots 21;t,x)$ completes the proof.
\hqedm

\bth
For any integers $k\ge0$, $\ell\ge1$, and $b\ge2$ we have
$$
S(1^k b^\ell;t,x)=
\frac{t^{k+\ell}x^{k+b\ell}}
{(1-x)^{k+1}\left((tx^b)^{\ell-1}(1-tx[b-1]_x
)+(1-x-tx)\dil\sum_{i=0}^{\ell-2}(1-x)^i(tx^b)^{\ell-2-i}\right)}.
$$
\eth
\pf
Suppose $w=w_1\ldots w_n\in\cS(1^k b^{\ell})$.  Then to have 
$1^k b^\ell$ as a suffix, we must have $w_n,\ldots,w_{n-\ell+1}\ge b$.

There are now two cases depending on the length of $w$.  If
$|w|=k+\ell$ then $w_1,\ldots,w_k$ are arbitrary positive integers.
If $|w|>k+\ell$ then write $w=yaz$ where $|z|=\ell$ and $a\in\bbP$.
In order to make sure that $1^k b^\ell$ does not have another
embedding intersecting $z$ it is necessary and sufficient that $a<b$.
And ruling out any embeddings inside $y$ is equivalent to 
$y\in\cA(1^k b^\ell)$.  We must also make sure that $|y|\ge k$ in
order to have $|w|>k+\ell$.

Let $S=S(1^k b^\ell;t,x)$ and $A=A(1^k b^\ell;t,x)$.
Turning all the information about $w$ into a generating function
identity gives
$$
S=\left(\frac{tx^b}{1-x}\right)^{\ell}
\left[ \left(\frac{tx}{1-x}\right)^k 
+ tx[b-1]_x \left(A -[k]_{tx/(1-x)}\right)
\right].
$$
Also, combining the two parts of Corollary~\ref{FSAtx} gives
$$
A=\frac{(1-x)(1-S)}{1-x-tx}.
$$
Substituting this expression for $A$ into our previous equation, one 
can easily solve for $S$ to obtain that 
\begin{equation*}
S(1^kb^{\ell};t,x)) = 
\frac{t^{k+\ell}x^{k+b\ell}(1-x-tx^b)}
{(1-x)^{k+1} ((1-x)^{\ell-1} (1-x-tx) + t^{\ell+1}x^{b\ell+1}[b-1]_x)}.
\end{equation*}
Thus to finish the proof, one need only show that
\begin{eqnarray*}
&&\frac{(1-x)^{\ell-1} (1-x-tx) + t^{\ell+1}x^{b\ell+1}[b-1]_x}
{(1-x-tx^b)}\\[5pt]
&&\hspace{10pt} = (tx^b)^{\ell-1}(1-tx[b-1]_x)
+(1-x-tx)\sum_{i=0}^{\ell-2}(1-x)^i(tx^b)^{\ell-2-i}
\end{eqnarray*}
which can be easily verified by cross multiplication.
\hqedm

\section{The M\"obius function}
\label{mf}

We will now show that the language for the
M\"obius function of ordinary factor order is not regular.
This is somewhat surprising because Bj\"orner and
Reutenauer~\cite{br:rmf} showed that this language is regular if one
considers ordinary subword order, and then Bj\"orner and
Sagan~\cite{bs:rmf} extended this result to generalized subword
order.  We will begin by reviewing some basic facts about M\"obius
functions.  The reader wishing more details can 
consult~\cite[Chapter 3]{sta:ec1}.

For any poset $P$, the {\it incidence algebra\/} of $P$ over the
integers is
$$
I(P)=\{\al:P\times P\ra\bbZ\ :\ \mbox{$\al(a,b)=0$ if $a\not\le b$}\}.
$$
This set is an algebra whose multiplication is given by {\it
  convolution\/}
$$
(\al*\be)(a,b)=\sum_{c\in P} \al(a,c)\be(c,b).
$$
It is easy to see that the identity for this operation is the
{\it Kronecker delta\/}
$$
\de(a,b)=\case{1}{if $a=b$,}{0}{else.}
$$
So it is possible for incidence algebra elements to have
multiplicative inverses.

One of the simplest elements of $I(P)$ is the {\it zeta function\/}
$$
\ze(a,b)=\case{1}{if $a\ge b$,}{0}{else.}
$$
Note that $F(u)$ can be rewritten as
$$
F(u)=\sum_{w\in P^*} \ze(u,w) w.
$$

It turns out that $\ze$ has a convolutional inverse $\mu$ in $I(P)$.
This function is important in enumerative and algebraic
combinatorics.  Bj\"orner~\cite{bjo:mff} has given a formula for $\mu$
in ordinary factor order which we will need.  To describe this result,
we must make some definitions.  
The {\it dominant outer factor\/} of $w$, denoted
$o(w)$, is the longest word other than $w$ which is both a prefix
 and a suffix of $w$.  Note that we may have $o(w)=\ep$.
The {\it dominant inner factor\/}
of $w=w_1\ldots w_\ell$, written $i(w)$, is $w_2\ldots w_{\ell-1}$.
Finally, a word is {\it flat\/} if all its elements are equal.  For
example, $w=abbaabb$ has $o(w)=abb$ and $i(w)=bbaab$.

\bth[Bj\"orner]
\label{factor}
In (ordinary) factor order, if $u\le w$ then
$$
\mu(u,w)=
\left\{\barr{ll}\mu(u,o(w))
       &\mbox{if $|w|-|u|>2$ and $u\le o(w)\not\le i(w)$,}\\
       1&\mbox{if $|w|-|u|=2$, $w$ is not flat, and $u=o(w)$ or $i(w)$,}\\
       (-1)^{|w|-|u|}&\mbox{if $|w|-|u|<2$,}\\
       0&\mbox{otherwise.\hspace*{245pt}\qed}\\
       \earr\right.
$$
\eth

Continuing the example
$$
\mu(b,abbaabb)=\mu(b,abb)=1.
$$
Note that this description is inductive.  It also implies that
$\mu(u,w)$ is $\pm1$ or $0$ for all $u,w$ in factor order.

We will show that the language
$\cM(u)=\{w\ :\ \mu(u,w)\neq0\}$ need not be regular.  To do this, we
will need the Pumping Lemma which we now state.  A proof can be found
in~\cite[pp.\ 55--56]{hu:iat}.
\ble[Pumping Lemma]
Let $\cL$ be a regular language.  Then there is a constant $n\ge1$ such
that any $z\in \cL$ can be written as $z=uvw$ satisfying
\ben
\item $|uv|\le n$ and $|v|\ge1$,
\item $uv^i w\in \cL$ for all $i\ge 0$.\hqed
\een
\ele
Roughly speaking, any word in a regular language has a  prefix of
bounded length such that pumping up the end of the prefix keeps one in
the language.

\bth
\label{cM}
Consider (ordinary) factor order where $P=\{a,b\}$.  Then $\cM(a)$ is
not regular.
\eth
\pf
Suppose, to the contrary, that $\cM(a)$ is regular and let $n$ be the
constant guaranteed by the pumping lemma.  We will derive a
contradiction by letting $z=a b^n a b^n a$ where, as usual, $b^n$
represents the 
letter $b$ repeated $n$ times.  

First we show that $z\in\cM(a)$.  Indeed, $o(z)=a b^n a$ and 
$i(z)= b^n a b^n$ which implies that $a\le o(z) \not\le i(z)$.  So we are in the
first case of Bj\"orner's formula and $\mu(a,z)=\mu(a,a b^n a)$.
Repeating this analysis with $a b^n a$ in place of $z$ gives
$\mu(a,z)=\mu(a,a)=1$.  Hence $z\in\cM(a)$ as promised.

Now pick any prefix $uv$ of $z$ as in the Pumping Lemma.  There are
two cases.  The first is if $u\neq\ep$.  
So $v=b^j$ for some $j$ with $1\le j<n$.  Picking $i=2$, we
conclude that $z'=uv^2w=a b^{n+j} a b^n a$ is in $\cM(a)$.
But $o(z')=a$ and 
$i(z')= b^{n+j} a b^n$.  Thus $|z'|-|a|> 2$ and $a\le o(z') \le i(z')$,
so $z'$ does not fall into any of the first three cases of Bj\"orner's
formula.  This implies that  $\mu(a,z')=0$ and hence $z'\not\in\cM(a)$, which is a
contradiction in this case.

The second possibility is that $u=\ep$ and $v=ab^j$ for some 
$0\le j<n$. Similar considerations to those in the previous paragraph
show that if we take $z'=uv^2w$ then $\mu(a,z')=0$ again.  So we have
a contradiction as before and the theorem is proved.\hqedm

\section{Comments, conjectures, and open questions}
\label{cco}

\subsection{Mixing factors and subwords}

It is possible to create languages using combinations of factors and subwords.
This is an idea that was first studied by Babson and
Steingr\'{\i}msson~\cite{bs:gpp} in the context of pattern avoidance
in permutations.  Many of the results we have proved can be
generalized in this way.  We will indicate how this can be done for
Theorem~\ref{NFA}. 

A {\it pattern\/} $p$ over P is a word in $P^*$ where certain pairs of
adjacent elements have been overlined (barred).  For example, in the pattern
$p=1\ol{133}24\ol{61}$ the pairs $13$, $33$, and $61$ have been overlined.
If $w\in P^*$ we will write $\ol{w}$ for the pattern where every pair
of adjacent elements in $w$ is overlined.  So every pattern has a
unique factorization of the form $p=\ol{y_1}\ \ol{y_2}\ \ldots\ \ol{y_k}$.  In the
preceding example, the factors are $y_1=1$, $y_2=133$, $y_3=2$,
$y_4=4$, and $y_5=61$. 

If $p=\ol{y_1}\ \ol{y_2}\ \ldots\ \ol{y_k}$ 
is a pattern and $w\in P^*$ then $p$
{\it embeds\/} into $w$, written $p\ra w$, if there is a subword
$w'=z_1 z_2\ldots z_k$ of $w$ where, for all $i$,
\ben
\item $z_i$ is a factor of $w$ with $|z_i|=|y_i|$, and
\item $y_i\le z_i$ in generalized factor order.
\een
For example $\ol{32}4\ra 14235$ and there is only one embedding,
namely $425$.  For any pattern $p$, define the language
$$
\cF(p)=\{w\in P^*\ :\ p\ra w\}
$$
and similarly for $\cS(p)$ and $\cA(p)$.  The next result generalizes
Theorem~\ref{NFA} to an arbitrary  pattern.
\bth
Let $P$ be any poset and let $p$ be a pattern over $P$.  Then there
are NFAs accepting $\cF(p)$, $\cS(p)$, and $\cA(p)$.
\eth
\pf
As before, it suffices to build an NFA, $\De$, for $\cS(p)$.  
It will be simplest to construct an NFA with $\ep$-moves, i.e., with
certain arcs labeled $\ep$ whose traversal does not append anything to
the word being constructed.  It is well known that the set of
languages accepted by NFAs with $\ep$-moves is still the set of
regular languages.

Let $p=\ol{y_1}\ \ol{y_2}\ \ldots\ \ol{y_k}$ 
be the  factorization of $p$ and, for all
$i$, let $\De_i$ be the automaton constructed in Theorem~\ref{NFA} for
$\cS(y_i)$.  We can paste these automata together to get $\De$ as
follows.   For each $i$ with $1\le i<k$, add an $\ep$-arc from every
final state of $\De_i$ to the initial state of $\De_{i+1}$.  Now let
the initial state of $\De$ be the initial state of $\De_1$ and the
final states of $\De$ be the final states of $\De_k$.
It is easy to see that the resulting NFA accepts the language $\cS(p)$.
\hqedm

\subsection{Rationality for infinite posets}

It would be nice to have a criterion that would imply rationality even
for some infinite posets $P$.  To this end, let
$\bx=\{x_1,\ldots,x_m\}$ be a set of commuting variables and consider the  formal power
series algebra $\bbZ[[\bx]]$.  Suppose we are given a 
function
$$
\wt:P\ra \bbZ[[\bx]]
$$
which then defines a weighting of words $w=w_1\ldots w_\ell\in P^*$ by
$$
\wt(w)=\prod_{i=1}^m \wt(w_i).
$$
To make sure our summations will be defined in $\bbZ[[\bx]]$, we
assume that there are only finitely many $w$ of any given weight and
call such a weight function {\it regular\/}.

For $u\in P^*$, let
$$
F(u;\bx)=\sum_{w\ge u} \wt(w)
$$
and similarly for $S(u;\bx)$ and $A(u;\bx)$.
Suppose we want to make sure that $S(u;\bx)$ is rational.  
As done in Section~\ref{pi}, we can consider a transfer matrix
with entries
$$
M_{T,U}=\sum_a \wt(a)
$$
where the sum is over all $a\in P$ occurring on arcs from $T$ to $U$.
Equation~\ree{LM} remains the same, so it suffices to make sure that
$M_{T,U}$ is always rational.

If there is an arc labeled $a$ from $T$ to $U$ then we must
have $U\sbe T'$ where $T'$ is given in equation~\ree{T'}.  Recalling
the definition of $\De$ from the proof of Theorem~\ref{NFA}, we see
that the $a$'s appearing in the previous sum are exactly those
satisfying
\ben
\item $a\ge u_{t+1}$ for $t+1\in U$, and
\item $a\not\ge u_{t+1}$ for $t+1\in T'-U$.
\een
To state these criteria succinctly, for any subword $y$ of $u$ we write $a\ge y$ 
(respectively, $a\not\ge y$) if 
$a\ge b$ (respectively, $a\not\ge b$) for all $b\in y$.  Finally, note that, from
the proof of Theorem~\ref{NFA}, similar transfer matrices can be
constructed for $F(u;\bx)$ and $A(u;\bx)$.  We have proved the
following result which generalizes Theorem~\ref{FSAtxthm}.
\bth
Let $P$ be a poset with a regular weight function $\wt:P^*\ra\bbZ[[\bx]]$,
and let $u\in P^*$.  Suppose that for any two
subwords $y$ and $z$ of $u$ we have
$$
\sum_{ \scriptstyle a\ge y \atop \rule{0pt}{8pt} \scriptstyle a\not\ge z} \wt(a)
$$
is a rational function.  Then so are $F(u;\bx)$, $S(u;\bx)$, and $A(u;\bx)$.\hqed
\eth

\subsection{Irrationality for infinite posets}

When $P$ is countably infinite it is possible for the generating functions we
have considered to be irrational.  As an example, pick a
distinguished element $a\in P$.  For any $A\sbe P$
with $a\in A$, we define an order $\le_A$ by insisting that the
elements of $P-\{a\}$ form an antichain, and that $a\le_A b$ if and
only if $b\in A$.  Consider the corresponding language $\cS_A$.
Clearly $\cS_A=(P-A)^* A$ and so no two of these languages are
equal.  It follows that the mapping $A\ra \cS_A$ is injective.
So one of the $\cS_A$ must be irrational since there are uncountably
many possible $A$ but only countably many rational functions in
$\bbZ\spn{\spn{P}}$. 

\subsection{Wilf equivalence and strong equivalence}

There are a number of open problems and questions raised by our work
on Wilf equivalence.

\medskip

(1)  {\bf  If $u\sim v$, then must $v$ be a rearrangement of $u$?}
     This is the case for all the Wilf equivalences we have proved.
     Note that if the answer is ``yes,'' then the Wilf equivalences
     for the symmetric groups
     given in Table 1 of Section~\ref{t} are actually Wilf
     equivalence classes.

\medskip

(2) {\bf What about Wilf equivalence in $[m]^*$ where $[m]=\{1,2,\ldots,m\}$?}
Given a positive integer $m$, one can define Wilf equivalence of words
$u,v\in[m]^*$ in the same way that we did for $\bbP^*$.  We write
$u\sim_m v$ for this relation.  Is it true that $u\sim_m v$ if and only
if $u\sim v$? 

\medskip

(3)  {\bf  If $u^+\sim v^+$ then is $u\sim v$?}  In other words, does
     the converse of Lemma~\ref{sim} (c) hold?  We note that the
     converse of (b) is true.  For suppose $1u\sim 1v$ and let
     $f:\cS(1u)\ra\cS(1v)$ be a corresponding map.  Then to construct
     $g:\cS(u)\ra\cS(v)$ we consider two cases for $w\in\cS(u)$.  If
     $|w|>|u|$ then 
     $w\in\cS(1u)$ so let $g(w)=f(w)$.  Otherwise $|w|=|u|$ and so
     let $g(w)=v+(w-u)$ where addition and subtraction is done
     componentwise.  It is easy to check that $g$ is well defined and
     weight preserving.

\medskip

(4)  {\bf  Find a theorem which, together with the results already
       proved, explains all the Wilf equivalences in $\fS_5$.}  
In particular, the results of Section~\ref{we} and the last
paragraph of Section~\ref{swe} generate all of the Wilf
equivalences in Table 1 with one exception.  In particular, our results show that
$$
31425 \sim 31524 \sim 42513 \sim 52413
\qmq{and}
32415 \sim 32514 \sim 41523 \sim 51423.
$$
but not why a permutation of the first group is Wilf equivalent to one
of the second.  However, we do have a conjecture which has been
verified by computer in a large number of examples and which would connect
these two groups.
\bcon
For any $a,b,c\in[2,\infty)$ we have
$$
a 1 b 2 c\sim a 2 b 1 c.
$$
\econ

\medskip

(5) {\bf  Is it always the case that the number of elements of $\fS_n$
    Wilf equivalent to a given permutation is a power of 2?}  This is
    always true in Table 1.

\medskip

(6) {\bf Is it true that $312\sim_s 213$?}  From our results on strong
    Wilf equivalence it follows that $12\sim_s 21$ and 
    $123 \sim_s 132 \sim_s 231 \sim_s 321$.  So all the Wilf
    equivalent elements in $\fS_2$ and $\fS_3$ are actually strongly
    Wilf equivalent with the possible exception of the pair in the
    question.  Of course, this breaks down in $\fS_4$ as noted in
    Section~\ref{swe}.

\medskip

(7)  {\bf  Does Theorem~\ref{rep} remain true if one replaces strong
     Wilf equivalence with ordinary Wilf equivalence throughout?}  If
     so, a completely different proof will have to be found 
     for that case.

\subsection{The language $\cM(u)$}

We have shown that $\cM(u)$ is not always regular and so the
corresponding generating function $M(u)$ is not always rational.  But this
leaves open whether $\cM(u)$  might fall into a more general class of
languages such as context free grammars.  A 
{\it context free  grammar\/}  or {\it CFG\/} is a quadruple
$G=(V,S,T,P)$ where 
\ben
\item $V$ is a finite set of {\it variables\/}, 
\item $S$ is a
special variable called the {\it start symbol\/}, 
\item $T$ is a finite set
of terminals disjoint from $V$, and 
\item $P$ is a finite set of 
{\it  productions\/}  of the form $A\ra\al$ where $A\in V$ and
$\al\in(V\cup T)^*$.
\een

There is a Pumping Lemma for CFGs, see~\cite[Section 6.1]{hu:iat}.  So
it is tempting to try and modify the proof of Theorem~\ref{cM} to show
that  $\cM(u)$ is not even a CFG.  However, all our attempts in that
direction have failed.  
Is $\cM(u)$ a CFG or not?

\section{Tables}
\label{t}

The following two tables were constructed by having a computer
calculate, for each composition $u$, the generating functions
$S(u;t,x)$.  This was done with the aid of the corresponding automaton
from Section~\ref{ca}.

\begin{center}
\renewcommand{\arraystretch}{1.25}
\begin{tabular}{|c|}
  \hline
  12, 21\\ \hline\hline
  123, 132, 231, 321\\ \hline
  213, 312\\ \hline\hline
  1234, 1243, 1342, 1432, 2341, 2431, 3421, 4321\\  \hline
  1324, 1423, 3241, 4231\\ \hline
  2134, 2143, 3412, 4312 \\  \hline
  3124, 3214, 4123, 4213\\  \hline
  2314, 2413, 3142, 4132\\ \hline\hline
  12345, 12354, 12453, 12543, 13452, 13542, 14532, 15432, \\23451, 23541, 24531,
25431, 34521, 35421, 45321, 54321\\  \hline
  12435, 12534, 14352, 15342, 24351, 25341, 43521, 53421 \\  \hline
  13245, 13254, 14523, 15423, 32451, 32541, 45231, 54231 \\  \hline
  21345, 21354, 21453, 21543, 34512, 35412, 45312, 54312 \\  \hline
  23145, 23154, 45132, 54132 \\ \hline
  32145, 32154, 45123, 54123\\ \hline
  24153, 25143, 34152, 35142\\ \hline
  14235, 14325, 15234, 15324, 42351, 43251, 52341, 53241\\ \hline
  31425, 31524, 32415, 32514, 41523, 42513, 51423, 52413\\ \hline
  24315, 25314, 41352, 51342\\ \hline
  24135, 25134, 43152, 53142\\ \hline
  34215, 35214, 41253, 51243\\ \hline
  34125, 35124, 42153, 52143\\ \hline
  41325, 42315, 51324, 52314\\ \hline
  41235, 43215, 51234, 53214\\ \hline
  42135, 43125, 52134, 53124\\ \hline
  13425, 13524, 14253, 15243, 34251, 35241, 42531, 52431\\ \hline
  21435, 21534, 43512, 53412\\ \hline
  24513, 25413, 31452, 31542\\ \hline
  23415, 23514, 41532, 51432\\ \hline
  31245, 31254, 45213, 54213\\  \hline\hline
  Table 1:  Wilf equivalences for permutations of at most 5 elements\\ \hline
\end{tabular}
\end{center}

\begin{center}
\renewcommand{\arraystretch}{1.75}
\begin{tabular}{|c|c|}
  \hline
  Equivalences & $S(u;t,x)$\\ \hline\hline
  1 & $\frac{tx}{1-x}$ \\ \hline
  2 & $\frac{tx^2}{(1-x)(1-tx)} $\\ \hline
  3 & $\frac{tx^3}{(1-x-tx+tx^3)} $\\ \hline
  11 & $\frac{t^2x^2}{(1-x)^2} $\\ \hline
  12,21 & $\frac{t^2x^3}{(1-x)^2(1-tx)} $\\ \hline
  13,31 & $\frac{t^2x^4}{(1-x)^2(1-tx-tx^2)}$\\ \hline
  22 & $\frac{t^2x^4}{(1-x)(1-x-tx+tx^2-t^2x^3)}$\\ \hline
  23,32 & $\frac{t^2x^5}{(1-x)(1-x-tx+tx^3-t^2x^4)}$\\ \hline
  33 & $\frac{t^2x^6}{(1-x)(1-x-tx+tx^3-t^2x^4-t^2x^5)}$\\ \hline
  111 & $\frac{t^3x^3}{(1-x)^3} $\\ \hline
  112,121,211 & $\frac{t^3x^4}{(1-x)^3(1-tx)} $\\ \hline
  122,221 & $\frac{t^3x^5}{(1-x)^2(1-x-tx+tx^2-t^2x^3)} $\\ \hline
  212 & $\frac{t^3x^5(1+tx^2)}{(1-x)(1-x+t^2x^3)(1-x-tx+tx^2-t^2x^3)} $\\ \hline
  113,131,311 & $\frac{t^3x^5}{(1-x)^3(1-tx-tx^2)} $\\ \hline
  213,312 & $\frac{t^3x^6(1+tx^3)}{(1-x)(1-x+t^2x^4)(1-x-tx+tx^3-t^2x^4)} $\\ \hline
  123,132,231,321 & $\frac{t^3x^6}{(1-x)^2(1-x-tx+tx^3-t^2x^4)} $\\ \hline
  222 & $\frac{t^3x^6}{(1-x)(1-2x-tx+x^2+2tx^2-tx^3-t^2x^3+t^2x^4-t^3x^5)} $\\ \hline
  133,331 & $\frac{t^3x^7}{(1-x)^2(1-x-tx+tx^3-t^2x^4-t^2x^5)} $\\ \hline
  313 & $\frac{t^3x^7(1+tx^3+tx^4)}{(1-x)(1-x+t^2x^4+t^2x^5)(1-x-tx+tx^3-t^2x^4-t^2x^5)} $\\ \hline
  223,232,322 & $\frac{t^3x^7}{(1-x)(1-2x-tx+x^2+tx^2+tx^3-tx^4-t^2x^4+t^2x^5-t^3x^6)} $\\ \hline
  323 & $\frac{t^3x^8(1+tx^3)}{(1-x)(1-2x-tx+x^2+tx^2+tx^3-tx^4-t^2x^4+t^2x^5-t^3x^6-t^3x^7+t^3x^8-t^4x^9-t^4x^{10})} $\\ \hline
  233,332 & $\frac{t^3x^8}{(1-x)(1-2x-tx+x^2+tx^2+tx^3-tx^4-t^2x^4+t^2x^6-t^3x^7)} $\\ \hline
  333 & $\frac{t^3x^9}{(1-x)(1-2x-tx+x^2+tx^2+tx^3-tx^4-t^2x^4+t^2x^6-t^3x^7-t^3x^8)} $\\  \hline\hline
\multicolumn{2}{|c|}{Table 2:  Wilf equivalences for $u$ with $|u|\le3$
  and $u_i\le3$ for all $i$.}
\\ \hline
\end{tabular}
\end{center}

\bigskip
\bibliographystyle{acm}
\begin{small}
\bibliography{ref}
\end{small}

\end{document}